\documentclass[10pt,letterpaper]{article}

\usepackage{times}
\usepackage{epsfig}
\usepackage{graphicx}
\usepackage{amsmath}
\usepackage{amssymb}

\usepackage{algorithm}
\usepackage{algorithmic}

\usepackage[margin=1in]{geometry}

\newcommand{\Reffig}[1]{Figure~\ref{fig:#1}}

\newcommand{\refsec}[1]{Sec.~\ref{sec:#1}}
\newcommand{\refeq}[1]{Eq.~\ref{eq:#1}}
\newcommand{\reffig}[1]{Fig.~\ref{fig:#1}}

\newcommand{\refalg}[1]{Alg.~\ref{alg:#1}}

\newcommand{\Frechet}[0]{Fr\'{e}chet }

\newcommand{\Exp}[0]{\operatorname{Exp}}
\newcommand{\Log}[0]{\operatorname{Log}}
\newcommand{\g}[0]{\mathfrak{g}}

\newcommand{\var}[0]{\operatornamewithlimits{Var}}

\renewcommand{\H}[0]{\mathcal{H}}

\newcommand{\ad}[0]{\operatorname{ad}}

\newcommand{\SO}[0]{\operatorname{SO}}
\newcommand{\so}[0]{\mathfrak{so}}

\newcommand{\etal}[0]{et~al. }

\usepackage[pagebackref=true,breaklinks=true,letterpaper=true,colorlinks,bookmarks=false]{hyperref}

\begin{document}

\title{Polynomial Regression on Riemannian Manifolds}

\author{Jacob Hinkle\and Prasanna Muralidharan \and P. Thomas Fletcher \and Sarang Joshi\\
Department of Bioengineering, University of Utah\\
72 Central Campus Drive, Salt Lake City, UT 84112\\
{\tt\small jacob@sci.utah.edu}
}

\maketitle

\begin{abstract}
    In this paper we develop the theory of parametric polynomial regression in
    Riemannian manifolds and Lie groups.  We show application of Riemannian
    polynomial regression to shape analysis in Kendall shape
    space.  Results are presented, showing the power of polynomial regression
    on the classic rat skull growth data of Bookstein as well as the analysis of the
    shape changes associated with aging of the corpus callosum from the OASIS Alzheimer's study.
\end{abstract}

\section{Introduction}

The study of the relationship between measured data and descriptive variables is
known as the field of regression analysis. As with most statistical techniques, regression analysis can be broadly divided into two classes: parametric and non-parametric.   
The most widely known parametric regression
methods are linear and polynomial regression in Euclidean space, 
wherein a linear or polynomial function is fit
in a least-squares fashion to observed data.  Such methods are the staple 
of modern data analysis. However,
classical regression suffers from the fundamental limitation that the
data must lie in a Euclidian vector space. 
The most common non-parametric regression approaches are kernel-based methods
and spline smoothing approaches which provide much more flexibility in the class
of regression functions.  However, their
non-parametric nature presents a challenge to inference problems, for instance
if one wishes to perform a hypothesis test to determine whether the trend for
one group of data is significantly different from that of another group. 

Fundamental to the analysis of anatomical imaging data within the framework of computational anatomy is the analysis of transformations and shape which are best represented as elements of Remaninan Manifolds rather than Euclidian vector spaces. In previous work, non-parametric kernel-based or spline-based methods have been
extended to observations that lie on a Remannian manifold with some
success~\cite{JDH:davis2010,JDH:jupp1987}, but intrinsic parametric
regression on Riemannian manifolds has received limited attention.  Most
recently,
Fletcher~\cite{JDH:fletcher2011} and
Niethammer~\etal~\cite{JDH:niethammer2011} have each independently
developed geodesic regression which generalizes
the notion of linear regression to Riemannian manifolds. 




The goal of the current work is to extend such work in order to accommodate more
flexibility in the model while remaining in the parametric setting.  The
increased flexibility introduced by the methods in this manuscript allow a
better description of the variability in the data, and ultimately will allow
more powerful statistical inference.  Our work builds off that of
Jupp~\&~Kent~\cite{JDH:jupp1987}, whose method for fitting parametric curves to
the sphere involved intricate unwrapping and rolling processes.  The work
presented in this paper allows one to fit regression curves on a general
Riemannian manifold, using intrinsic methods and avoiding the need for
unwrapping and unrolling.

We demonstrate the usefulness of our algorithm in three studies of shape.  By
applying our algorithm to Bookstein's classical rat skull growth
dataset~\cite{JDH:bookstein1991}, we show that we are able to obtain a
parametric regression curve of similar quality to that produced by
non-parametric methods~\cite{JDH:kent2001}.  We also demonstrate in a 2D corpus
callosum aging study that, in addition to providing more flexibility in the
traced path, our polynomial model provides information about the optimal
parametrization of the time variable.



\section{Methods}

\subsection{Preliminaries}

Let $(M,g)$ be a Riemannian manifold~\cite{JDH:docarmo1992}.  For each point $p\in M$, the metric $g$
determines an inner product on the tangent space $T_pM$ as well as a way to
differentiate vector fields $X,Y$ with respect to one another.  That derivative
is referred to as the covariant derivative and is denoted $\nabla_{X} Y$.  If
$V\in\mathfrak{X}(M)=\{f\in C^\infty(M,TM):f(p)\in T_pM,\forall p\in M\}$ is a
smooth vector field on $M$ and $\gamma:[0,T]\to M$ is a smooth curve on $M$ then
the covariant derivative of $V$ along $\gamma$ is the time derivative of $V$ in
a reference frame along $\gamma$
\begin{align}
     \nabla_{\dot{\gamma}(t)} V &= \frac{D}{dt}V(\gamma(t)),
\end{align}
where the intrinsic derivative $\frac{D}{dt}$ is determined by the metric. 
Geodesics $\gamma:[0,T]\to M$ are characterized by the second-order
covariant differential equation (c.f. \cite{JDH:docarmo1992})
\begin{align}
    \nabla_{\dot{\gamma}} \dot{\gamma} &= 0.
\end{align}
This equation, called the geodesic equation, uniquely determines
geodesics up to a choice of initial conditions $(\gamma(0),\dot{\gamma}(0))\in TM$.
The mapping from the tangent bundle into the manifold is called the
exponential map, $\Exp:TM\to M$.  Fixing the base
point $p\in M$, the exponential map is injective on a zero-centered ball
$B$ in $T_pM$ of some non-zero (possibly infinite) radius.  Thus for a
point $q$ in the image of $B$ under $\Exp_p$ there exists a unique
vector $v\in T_pM$ corresponding to a minimal length path under the
exponential map from $p$ to $q$.  The mapping of such points $q$
to their associated tangent vectors $v$ at $p$ is called the log map of $q$ at
$p$, denoted $v = \Log_p q$.

Given a curve $\gamma:[0,T]\to M$ we'll want to relate tangent vectors
at different points along the curve.  These relations are governed
infinitesimally by the covariant derivative $\nabla_{\dot{\gamma}}$.  A vector
field $X:[0,T]\to T_{\gamma(t)}M$ along the curve $\gamma$ is parallel
transported along $\gamma$ if it satisfies the parallel transport
equation:
\begin{align}
    \nabla_{\dot{\gamma}} X(t) &= 0
\end{align}
for all times $t\in[0,T]$.  Notice that the geodesic equation is a
special case of parallel transport, in which we require that the
velocity is parallel transported along the curve itself.

\subsection{Riemannian Polynomials}

Given a vector field $X$ along a curve $\gamma$, the covariant
derivative of $X$ gives us a way to define vector fields which are
``constant'' along $\gamma$, as parallel transported vectors.  Geodesics
are generalizations to the Riemannian manifold setting of curves with
constant first derivative.

The covariant derivative of a vector field along a curve gives another
vector field along the curve.  We will apply that derivative repeatedly
to examine curves with constant higher derivatives.  For instance, we
refer to the vector field $\nabla_{\dot{\gamma}(t)}\dot{\gamma}(t)$ as
the acceleration of the curve $\gamma$.  Curves with constant
acceleration are generalizations of quadratic curves in $\mathbb{R}$ and
satisfy the second order polynomial equation
\begin{align}
    (\nabla_{\dot{\gamma}})^2\dot{\gamma}(t) &= 0 , \\
     \mbox{with initial conditions:~}\qquad
     \gamma(0)\ , \ \dot{\gamma}(0) & \mbox{ and } \ddot{\gamma}(0) . \nonumber
\end{align}
Extending this idea, a cubic polynomial is defined as curve having constant jerk
(time derivative of acceleration), and so on.  Generally, a $k$th order
polynomial in $M$ is defined as a curve $\gamma:[0,T]\to M$ satisfying
\begin{align}
    (\nabla_{\dot{\gamma}})^k\dot{\gamma}(t) &= 0 \\
    \mbox{with initial conditions:~}\qquad     \gamma(0) &\mbox{ and }
    \dot{\gamma}^{i}(0), \quad i=1,\cdots,k,  \nonumber
\end{align}
for all time points $t\in[0,T]$. As with polynomials in Euclidian space, a
$k^\mathrm{th}$ order polynomial is determined by the $(k+1)$ initial conditions at $t=0$.

The covariant differential equation governing the evolution of Riemannian
polynomials is linearized in the same way that a Euclidian ordinary differential
equation is.  Introducing vector fields $v_1(t),\dots,v_{k}(t)\in
T_{\gamma(t)}M$, we can write the system of covariant differential equations as
\begin{align}
    \dot{\gamma}(t) &= v_1(t) \\
    \nabla_{\dot{\gamma}}v_1(t) &= v_2(t) \\
    &\vdots \\
    \nabla_{\dot{\gamma}}v_{k-1}(t) &= v_{k}(t) \\
    \nabla_{\dot{\gamma}}v_{k}(t) &= 0.
\end{align}

\begin{algorithm}[ht]
    \caption{\label{alg:forward}Pseudocode for forward integration of
    $k^\mathrm{th}$ order Riemannian polynomial}
    \begin{algorithmic}
        \STATE $\gamma \gets \gamma(0)$
        \FOR {$i=1,\ldots,k$}
          \STATE $v_i \gets v_i(0)$
        \ENDFOR
        \STATE $t\gets 0$
        \REPEAT
            \STATE $w \gets v_1$
            \FOR {$i=1,\ldots,k-1$}
            \STATE $v_i \gets \operatorname{ParallelTransport}_\gamma(\Delta t w,v_i + \Delta t v_{i+1})$
            \ENDFOR
            \STATE $v_k \gets \operatorname{ParallelTransport}_\gamma(\Delta t w,v_k)$
            \STATE $\gamma \gets \Exp_\gamma(\Delta t w)$
            \STATE $t\gets t+\Delta t$
        \UNTIL {t=T}
    \end{algorithmic}
\end{algorithm}

\begin{figure}[th]
    \begin{center}
        \includegraphics[width=0.5\textwidth,trim=90 95 90 80,clip=true]{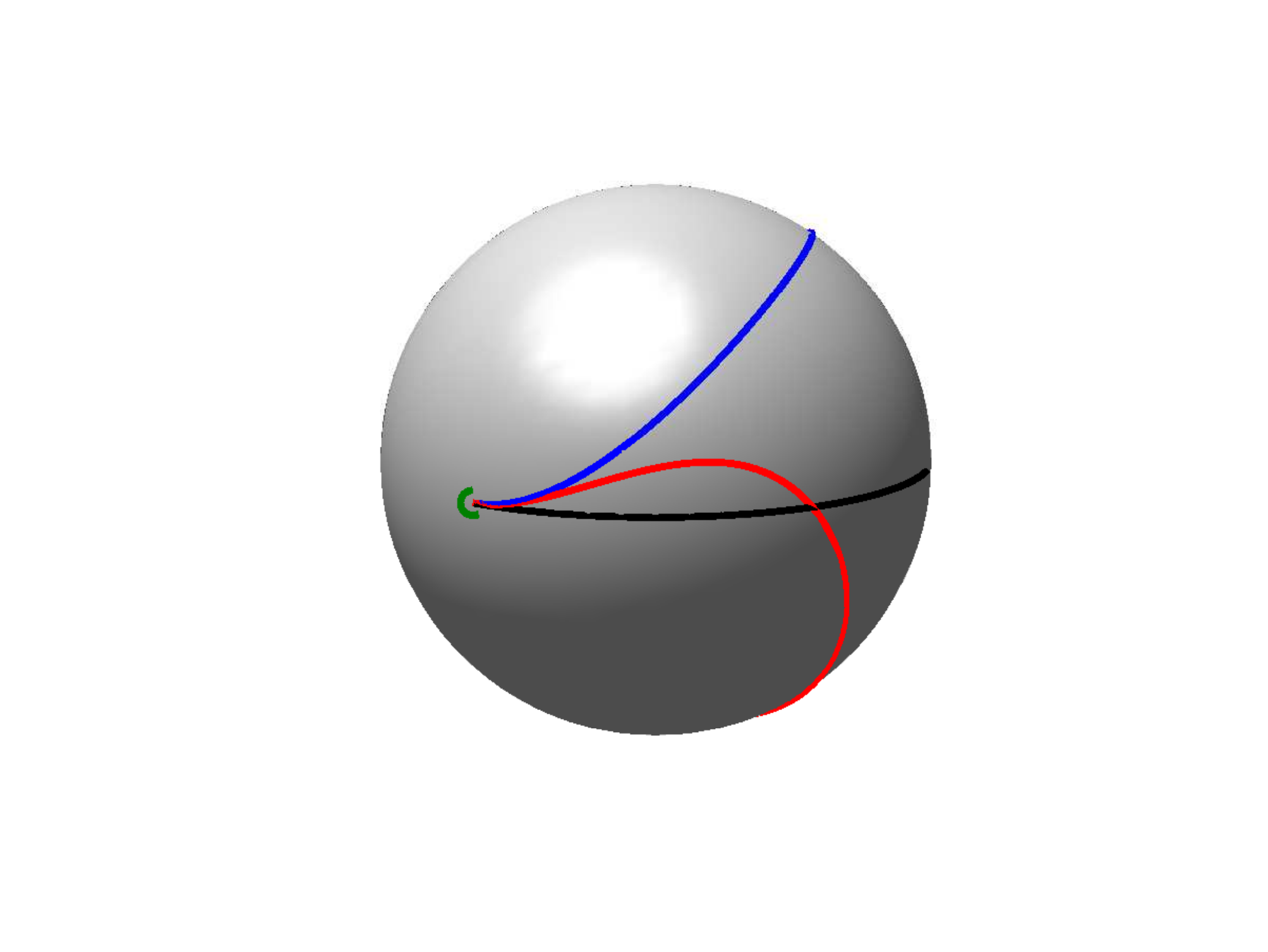}
    \end{center}
    \caption{Sample polynomial curves emanating from a common basepoint (green)
    on the sphere (black=geodesic,blue=quadratic,red=cubic).}
    \label{fig:spherecurves}
\end{figure}

The Riemannian polynomial equation cannot, in general, be solved in closed form,
and must be integrated numerically.  In order to discretize this system of
covariant differential equations, we implement the covariant integrator depicted
in Alg.~\ref{alg:forward}.  At each step of the integrator, each vector is
incremented within the tangent space at $\gamma(t)$ and the results are parallel
transported infinitesimally along a geodesic from $\gamma(t)$ to
$\gamma(t+\Delta t)$.  The only ingredients necessary to integrate a polynomial
forward in time are the exponential map and parallel transport on the manifold.

\Reffig{spherecurves} shows the result of integrating polynomials of order one,
two, and three.  Initial velocity, acceleration, and jerk were chosen and a
cubic polynomial integrated to obtain the blue curve.  Then the initial jerk was
set to zero and the blue quadratic curve was integrated, followed by the black
geodesic whose acceleration was set to zero.

\subsection{Estimation of Parameters for Regression}
\label{sec:estimation}

In order to regress polynomials against observed data $y_i\in M,i=1,\dots,N$,
we define the following constrained objective function
\begin{align}
    E_0(\gamma(0),v_1(0),\dots,v_k(0)) &= \frac{1}{N}\sum_{i=1}^N d(\gamma(t_i),y_i)^2 \\
    \mbox{subject to:}\qquad\qquad \dot{\gamma}(t) &= v_1(t) \\
    \nabla_{\dot{\gamma}}v_1(t) &= v_2(t) \\
    &\vdots \\
    \nabla_{\dot{\gamma}}v_{k-1}(t) &= v_{k}(t) \\
    \nabla_{\dot{\gamma}}v_{k}(t) &= 0.
\end{align}
which is minimized in order to find the optimal initial conditions
$\gamma(0),v_i(0),i=1,\dots,k$, which we will refer to as the parameters of our
model. 

In order to determine the optimal parameters $\gamma(0),v_i(0),i=1,\dots,k$, we introduce
Lagrange multiplier vector fields $\lambda_i\in \mathfrak{X}(M), i=0,\dots,k$,
often called the adjoint variables, and define the unconstrained objective
function
\begin{align}
    E(&\gamma(0),v(0),\lambda) = \frac{1}{N}\sum_{i=1}^N d(\gamma(t_i),y_i)^2 \\
    &+\int_0^T \langle\lambda_0(t), \dot{\gamma}(t)-v_1(t)\rangle dt \\
    &+\sum_{i=1}^{k-1}\int_0^T \langle\lambda_i(t), \nabla_{\dot{\gamma}}v_i(t) - v_{i+1}(t) \rangle dt \\
    &+\int_0^T\langle\lambda_k(t), \nabla_{\dot{\gamma}}v_k(t)\rangle
    dt.
    \label{eq:Eunconstrained}
\end{align}
As is standard practice, the optimality conditions for this equation are
obtained by taking variations with respect to all arguments of $E$, integrating
by parts when necessary.  The resulting variations with respect the the adjoint
variables yield the original dynamic constraints: the polynomial equations.
Variations with respect to the primal variables gives rise to the following
system of equations, termed the adjoint equations.  The adjoint
equations take the following form (see appendix for derivation):
\begin{align}
    \nabla_{\dot{\gamma}}\lambda_0  &= -\frac{2}{N}\sum_{i=1}^N
    \delta(t-t_i) \Log_{\gamma}y_i - \sum_{i=1}^kR(v_i,\lambda_i)v_1\\
    \nabla_{\dot{\gamma}}\lambda_1 &= -\lambda_0 \\
    &\vdots \\
    \nabla_{\dot{\gamma}}\lambda_{k} &= -\lambda_{k-1}
    \label{eq:adjointeqns},
\end{align}
where $R$ is the Riemannian curvature tensor and the Dirac $\delta$
functional indicates that the order zero adjoint variable takes on jump
discontinuities at time points where data is present.  Gradients of $E$
with respect to initial and final conditions give rise to the terminal
endpoint conditions for the adjoint variables, as well as expressions
for the gradients with respect to initial conditions.
\begin{align}
    \lambda_i(T) &= 0, i=0,\cdots,k \\
    \delta_{\gamma(0)} E &= -\lambda_0(0) \\
    \delta_{v_i(0)} E &= -\lambda_i(0)
\end{align}
In order to
determine the value of the adjoint vector fields at $t=0$, and thus the
gradients of our functional $E_0$, the adjoint variables are first
initialized to zero at time $T$, then the system \refeq{adjointeqns} is
evolved backward in time to $t=0$.

Given the gradients with respect to the initial conditions, a simple
steepest descent algorithm is used to optimize the functional.  The
update to $\gamma(0)$ is computed using the exponential map, and the
vectors $v_i(0)$ are updated via parallel translation.  

Note that in the special case of a zero-order polynomial ($k=0$), the
only gradient $\lambda_0$ is simply the mean of the log map vectors at
the current estimate of the \Frechet mean.  So this method generalizes
the common method of \Frechet averaging on
manifolds~\cite{JDH:fletcher2004}.

The curvature term, in the case $k=1$, indicates that $\lambda_1$ is a
sum of adjoint Jacobi fields.  So this approach subsumes geodesic
regression as presented by Fletcher~\cite{JDH:fletcher2011}.  For higher
order polynomials, the adjoint equations represent a generalized Jacobi
system. 

In practice, it is often not necessary to explicitly compute the
curvature terms.  In the case that the manifold $M$ embeds into a
Hilbert space, the extrinsic adjoint equations can be computed by taking
variations in the ambient space, using standard methods.  Such an
approach gives rise to the regression algorithm found in
Niethammer~\etal~\cite{JDH:niethammer2011}, for example.

\begin{algorithm}[ht]
    \caption{\label{alg:adjoint}Pseudocode for reverse integration of
    adjoint equations for $k^\mathrm{th}$ order Riemannian polynomial}
    \begin{algorithmic}
        \STATE $\gamma \gets \gamma(T)$
        \FOR {$i=0,\ldots,k$}
          \STATE $\lambda_i \gets 0$
        \ENDFOR
        \STATE $t\gets T$
        \REPEAT
            \STATE $w \gets v_1(t)$
            \STATE $\lambda_0 \gets \lambda_0+ \Delta t \sum_{i=1}^kR(v_i,\lambda_i)v_1$
            \IF{$t=t_i$}
                \STATE $\lambda_0 \gets\lambda_0 +\frac{2}{N}\Log_\gamma y_i$
            \ENDIF
            \FOR {$i=k,\ldots,1$}
            \STATE $\lambda_i \gets
            \operatorname{ParallelTransport}_\gamma(-\Delta t
            w,\lambda_i + \Delta t \lambda_{i-1})$
            \ENDFOR
            \STATE $\lambda_0 \gets
            \operatorname{ParallelTransport}_\gamma(-\Delta t w,\lambda_0)$
            \STATE $\gamma \gets \Exp_\gamma(-\Delta t w)$
            \STATE $t\gets t-\Delta t$
        \UNTIL {t=0}
        \STATE $\delta_{\gamma(0)}E \gets -\lambda_0$
        \FOR{$i=1,\ldots,k$}
            \STATE $\delta_{v_i(0)}E \gets -\lambda_i$
        \ENDFOR
    \end{algorithmic}
\end{algorithm}

\subsection{Time Reparametrization}
\label{sec:reparam}

In the geodesic model, curves propagate at a constant speed, a result of their
extremal action property.  However, using higher-order polynomials it is
possible to generate curves whose images match those of a geodesic, but whose
time-dependence has been reparametrized.  If the initial conditions consist of
all $v_i$ collinear, then this will necessarily be the case.  Polynomials
provide flexibility not only in the class of paths that are possible, but in the
time dependence of the curves traversing those paths.  Regression models could
even be imagined in which the operator wishes to estimate geodesic paths, but is
unsure of parametrization, and so enforces the estimated parameters to be
collinear.


\subsection{Polynomials on Riemannian Lie Groups}

A special case arises when the manifold is a Lie group and the metric is
left-invariant.  In this case, left-translation by $\gamma$ and $\gamma^{-1}$, denoted
$L_\gamma, L_{\gamma^{-1}}$, allow one 
to isometrically represent tangent vectors at $\gamma$ as tangent vectors at the identity,
using the pushforward $L_{\gamma^{-1}*}$.
The tangent space at the identity is isomorphically identified with the Lie algebra $\g$.  The algebraic operation in $\g$ is called
the Lie bracket $[,]:\g\times\g\to\g$.  For $X,Y\in\g$, this bracket operation,
also called the adjoint action of $X$ on $Y$ and often denoted $\ad_XY=[X,Y]$, is
bilinear and alternating in its arguments (i.e. $\ad_XY=-\ad_YX$).  Fixing $X$,
the operator $\ad_X$ is linear and its adjoint is computed with respect to
the metric:
\begin{align}
    \langle \ad_X Y, Z\rangle &= \langle Y, \ad_X^\dagger Z\rangle.
\end{align}
We will refer to the operator $\ad_X^\dagger$ as the adjoint-transpose action of
$X$.

Suppose $X$ is a vector field along the curve $\gamma$ and
$X_c=L_{\gamma^{-1}*}X\in\g$ is its representative at the identity.
Similarly, $\omega_c=L_{\gamma^{-1}*}\dot{\gamma}$ is the representative
of the curve's velocity in the Lie algebra.  Then the covariant derivative of
$X$ along $\gamma$ evolves in the Lie algebra via 
\begin{align}
    L_{\gamma^{-1}*}&\nabla_{\dot{\gamma}}X = \dot{X_c} - \frac{1}{2}\left( \ad_{\omega_c}^\dagger X_c +\ad_{X_c}^\dagger
    \omega_c - \ad_{\omega_c}X_c \right).\label{eq:algebracovderiv}
\end{align}
In the special case when $X=\dot{\gamma}$ so that $X_c=\omega_c$, if
we set this covariant derivative equal to zero we have the geodesic equation:
\begin{align}
    \dot{\omega_c} &= \ad_{\omega_c}^\dagger \omega_c,
\end{align}
which in this form is often refered to as the Euler-Poincare
equation~\cite{JDH:arnold1989,JDH:holm1998}.  The curvature tensor is given by the standard formula
\begin{align}
    R(X,Y)Z &= \nabla_X \nabla_Y Z - \nabla_Y\nabla_X Z - \nabla_{[X,Y]}
    Z.
\end{align}
Applying \refeq{algebracovderiv}, in a Lie group with left-invariant
metric, these terms are computed in the Lie algebra using the formulas:
\begin{align}
    \nabla_X\nabla_Y Z &= \frac{1}{4}\big(
    \ad_{\ad_X Y}Z - \ad_X\ad_Z^\dagger Y - \ad_X\ad_Y^\dagger Z \\
    &\qquad-\ad_{\ad_Y Z}^\dagger X + \ad_{\ad_Z^\dagger Y}^\dagger X + \ad_{\ad_Y^\dagger Z}^\dagger X\\
    &\qquad-\ad_X^\dagger{\ad_Y Z} + \ad_X^\dagger\ad_Z^\dagger Y + \ad_X^\dagger\ad_Y^\dagger Z
    \big) \\
    \nabla_{[X,Y]}Z &= \frac{1}{2}\left(\ad_{\ad_X Y}Z - \ad_{Z}^\dagger\ad_X Y -
    \ad_{\ad_X Y}^\dagger Z\right).\label{eq:algebracurvature}
\end{align}
Also note that if the metric is both left- and right-invariant, then the
adjoint-transpose action is alternating and we can simplify the
covariant derivative further.  In particular, in the presence of such a
bi-invariant metric, the first few covariant derivatives take the
following forms:
\begin{align}
    L_{\gamma^{-1}*}\left(\nabla_{\dot{\gamma}}\right)^k\dot{\gamma} &=
    \left( \frac{d}{dt}+\frac{1}{2}\ad_{\omega_c} \right)^k\omega_c \\
    L_{\gamma^{-1}*}\nabla_{\dot{\gamma}}\dot{\gamma} &=
    \dot{\omega_c} \\
    L_{\gamma^{-1}*}\left(\nabla_{\dot{\gamma}}\right)^2\dot{\gamma}&=
    \ddot{\omega_c} + \frac{1}{2}\ad_{\omega_c}\dot{\omega_c}
    \\
    L_{\gamma^{-1}*}\left(\nabla_{\dot{\gamma}}\right)^3\dot{\gamma}&=\dddot{\omega_c}+\ad_{\omega_c}\ddot{\omega_c}
    + \frac{1}{4}\ad_{\omega_c}\ad_{\omega_c}\dot{\omega_c}.
\end{align}
This allows us to solve the forward evolution in the convenient setting of the
vector space $\g$.  Parallel transport requires solution of the equation
\begin{align}
    \left( \frac{d}{dt}+\ad_{\omega_c} \right)X &= 0,
\end{align}
which can be integrated using Euler integration or similar time-stepping
algorithms, so long as the adjoint action can be easily computed.  
Because of the simplified covariant derivative formula for a
bi-invariant metric, the curvature takes the simple
form~\cite{JDH:docarmo1992}:
\begin{align}
    R(X,Y)Z&= \frac{1}{4}\ad_Z\ad_X Y
\end{align}

\subsection{Coefficient of Determination for Regression in Metric Spaces}

It will be useful to define a statistic which will indicate how well our model
fits some set of observed data.  As in \cite{JDH:fletcher2011}, we compute the
coefficient of determination of our regression polynomial $\gamma(t)$, denoted
$R^2$.  The first step to computing $R^2$ is to compute the variance of the
data.  As discussed in \cite{JDH:fletcher2011}, the natural choice of total
variance statistic is the \Frechet variance, defined by
\begin{align}
    \var\{y_i\} &= \frac{1}{N} \min_{\bar{y}\in M} \sum_{i=1}^N
    d(\bar{y},y_i)^2.
\end{align}
Note that the \Frechet mean $\bar{y}$ itself is the 0-th order polynomial
regression against the data $\{y_i\}$ and the variance is the value of the
objective function $E_0$ at that point.  We also define the sum of squared error
for a curve $\gamma$ as the value $E_0(\gamma)$:
\begin{align}
    SSE &= \frac{1}{N} \sum_{i=1}^N d(\gamma(t_i),y_i)^2
\end{align}
Then the coefficient of determination is defined as
\begin{align}
    R^2 &= 1- \frac{SSE}{\var\{y_i\}}
\end{align}
Any regressed polynomial will have $SSE<\var\{y_i\}$, $R^2$ will have
a value between 0 and 1, with 1 indicating a perfect fit and 0 indicating that
the curve $\gamma$ provides no better fit than does the \Frechet mean.


\section{Examples}

\subsection{The $n$-Dimensional Sphere}

Suppose $M=S^n=\{p\in \mathbb{R}^{n+1}:\|p\|=1\}$.  Then geodesics are great
circles and the exponential map is given by
\begin{align}
    \Exp_p v &= \cos\theta\cdot p +
    \sin\theta\frac{v}{\|v\|},\quad \theta = \|v\|
\end{align}
The corresponding log map is the inverse of this function:
\begin{align}
    \Log_p q &=
    \theta\frac{q-(p^Tq)p}{\|q-(p^T q)p\|},\quad\theta = \cos^{-1}(p^Tq).
\end{align}
The Riemannian curvature tensor is \cite{JDH:docarmo1992}
\begin{align}
    R(X,Y)Z &= (X^TZ)Y -(Y^TZ)X.
\end{align}
Parallel translation of a vector $X$ along the exponential map of a vector $v$
is performed as follows.  The vector $X$ is decomposed into a vector parallel to
$v$, which we denote $X^v$, and a vector orthogonal to $v$, $X^\perp$.
$X^\perp$ is left unchanged by parallel transport along $\Exp_pv$, while
$X^v$ transforms in a similar way as $v$:
\begin{align}
    X^v\mapsto \left(-\sin\theta p  + \cos\theta
    \frac{v}{\|v\|}\right)\frac{v^T}{\|v\|} X^{v}.
\end{align}
Using this exponential map and parallel transport, the integrator depicted in
Alg.~\ref{alg:forward} was implemented.  \Reffig{spherecurves} shows the result
of integrating polynomials of order one, two, and three, using the equations
from this section, for the special case of a two-dimensional sphere.

\subsection{The Lie Group $SO(3)$}
In order to illustrate the algorithm in a Lie group, we consider the Lie group
$\SO(3)$ of orthogonal 3-by-3 real matrices with determinant one.  We review the
basics of $\SO(3)$ here, but the reader is encouraged to consult
\cite{JDH:arnold1989} for a full treatment of $\SO(3)$ and derivation of the
Euler-Poincare equation there.

The Lie algebra for the rotation group is $\so(3)$, the set of all
skew-symmetric 3-by-3 real matrices.  Such matrices can be identified with
vectors in $\mathbb{R}^3$ by the identification
\begin{align}
    x=\left( \begin{array}[c]{c}
        a \\
        b \\
        c
    \end{array}\right) &\mapsto \sigma(x) =\left( \begin{array}[c]{ccc}
        0 & -c & b \\
        c & 0 & -a \\
        -b & a & 0
    \end{array}\right)
\end{align}
under which the Lie bracket takes the convenient form of the cross-product of
vectors:
\begin{align}
    x\times y = \sigma(x)y &\mapsto \sigma(x\times y)
    = [ \sigma(x), \sigma(y) ].
\end{align}
With the $\mathbb{R}^3$ representation of vectors in $\so(3)$, an invariant
metric is determined by the choice of a symmetric, positive-definite matrix $A$:
\begin{align}
    \langle x,y\rangle &= x^T A y, \forall x,y\in \so(3)\cong\mathbb{R}^3.
\end{align}
With this metric the adjoint-transpose action can be written as
\begin{align}
    \ad_x^\dagger y &= -A^{-1} \left(x\times Ay\right).
\end{align}
If we plug this into the general Lie group formula we obtain the geodesic
equation in the Lie algebra:
\begin{align}
    \dot{\omega_c} &= -A^{-1}(\omega_c\times A\omega_c)
\end{align}
which must be integrated numerically to find $\omega_c$ at each time.  To
compute the exponential map $\Exp_pv$, the algorithm depicted in
\refalg{forward} is employed.  The matrix $v$ is converted to a skew-symmetric
matrix by multiplication by $p^{-1}$, which translates $v$ back into the Lie
algebra.  The resulting vector is converted to a vector $\omega_c(0)$.  The
geodesic equation is then integrated to find $\omega_c$ at future times.  The
point $p$ then evolves along the exponential map by time-stepping, where at each
step $\omega_c$ is converted back into a skew-symmetric matrix $W$, multiplied
by $\Delta t$ and exponentiated to obtain a rotation matrix.  The Euler step is
then given by left-multiplication by the resulting matrix.

Plugging into \refeq{algebracurvature}, the terms need to compute the curvature
of $\SO(3)$ are
\begin{align}
    \nabla_X\nabla_Y Z &= \frac{1}{4}\big(
    (X\times Y)\times Z + X\times A^{-1}(Z\times AY) + X\times
    A^{-1}(Y\times AZ)\\
    &\qquad+A^{-1}( (Y\times Z)\times AX)+ A^{-1}( A^{-1}(Z\times AY) \times AX )+ A^{-1}( A^{-1}(Y \times AZ)\times AX) \\
    &\qquad+A^{-1}( X\times A(Y\times Z))+ A^{-1}(X\times (Z\times AY))+ A^{-1}(X\times (Y\times AZ))
    \big) \\
    \nabla_{[X,Y]}Z &= \frac{1}{2}\left( (X\times Y)\times Z + A^{-1}(Z\times A(X\times Y)) +A^{-1}( (X\times Y)\times AZ)\right),
\end{align}
and the parallel transport equation is given by
\begin{align}
    \dot{X_c} &= \frac{1}{2}\left(A^{-1}\left(-X_c\times A\omega_c-
    \omega_c\times AX_c\right) - \omega_c\times X_c\right).
\end{align}
The evolution equations are then integrated numerically as in \refalg{forward}.

\subsection{Kendall Shape Space}

The analysis of shape is a common challenge in medical imaging.  Shape is 
defined as the property of the data that is invariant to scale or rotation.
It was in
this setting that Kendall~\cite{JDH:kendall1984} originally developed his theory
of shape space.  We describe here Kendall's shape space of $m$-landmark
point sets in $\mathbb{R}^d$, denoted $\Sigma_d^m$.  For a more complete overview
of Kendall's shape space, including computation of the curvature tensor, the reader is encouraged to consult the work of
Kendall and Le~\cite{JDH:kendall1984,JDH:kendall1989,JDH:le1993}.

Let $x=(x_i)_{i=1,\ldots,m},x_i\in\mathbb{R}^d$ be a labelled pointset in
$\mathbb{R}^d$.  Translate the pointset so that the centroid resides at 0,
removing any effects due to translations.  Another degree of freedom due to
scaling is removed by requiring that $\sum_{i=1}^m \|x_i\|^2=1$.
After this standardization, one is left with a point in $S^{md-1}$, commonly
referred to as \emph{preshape space} in this context.

Kendall shape space $\Sigma_d^m$ is obtained by taking the quotient of the
preshape space by the natural action of the rotation group $SO(d)$ on pointsets.
In practice, points in the quotient (refered to as \emph{shapes}) are
represented by a member of their equivalence class in preshape space.  The work
of O'Neill~\cite{JDH:oneill1966} characterizes the link in geometry between
the shape and preshape spaces.  We describe now how to compute exponential maps,
log maps, and parallel transport in shape space, using representatives in
$S^{md-1}$.


First, notice that there exist tangent vectors in $S^{md-1}$ corresponding to
global rotations of the preshape.  At a particular point $p$ in preshape space,
these vectors form a linear subspace of the tangent space $T_pS^{md-1}$, which
we define as the vertical subspace.  Curves moving along vertical tangent
vectors result in rotations of a preshape, and so do not indicate any change in
actual shape.  Since we are concerned with shape change and wish to ignore
effects due to rotation, it will be important to project such tangent vectors to
be purely non-vertical, or horizontal.  We denote the horizontal projection by
$\H:TS^{md-1}\to TS^{md-1}$.


A vertical vector in preshape space arises as the derivative of a rotation of a
preshape.  The derivative of such a rotation is a skew-symmetric matrix $W$, and
its action on a preshape $x$ has the form $(Wx_1,\ldots,Wx_n)\in TS^{md-1}$.
The vertical subspace is then spanned by such tangent vectors arising from any
linearly independent set of skew-symmetric matrices.  The projection $\H$ is
performed by taking such a spanning set, performing Gram-Schmidt
orthonormalization, and removing each component.

The horizontal projection allows us to relate the covariant derivative on the
sphere to that on shape space.  Lemma 1 of O'Neill~\cite{JDH:oneill1966} states
that if $X,Y$ are horizontal vector fields at some point $p$ in preshape
space, then
\begin{align}
    \H\nabla_XY &= \nabla_{X^*}^*Y^*,
\end{align}
where $\nabla$ denotes the covariant derivative on preshape space and
$\nabla^*,X*$, and $Y*$ are their counterparts in shape space.

For a general shape space $\Sigma_d^m$, the exponential map and parallel
translation are performed using representatives in $S^{md-1}$.  In practice,
this usually must be done in a time-stepping algorithm, in which at each time
step an infinitesimal spherical parallel transport is performed, followed by the
horizontal projection.  The resulting algorithm can be used to compute the
exponential map as well.  Computation of the log map is less trivial, as it
requires an iterative optimization routine.  The target shape is first aligned
to the base shape using Procrustes alignment.  Then, at each iteration of the
log map routine,  the exponential map is integrated forward, compared to the
target preshape, then the resulting tangent vector is adjoint 
parallel transported back to the base point, for update in a steepest descent iteration.

A special case arises in the case when $d=2$.  In this case the exponential map,
parallel transport and log map have closed form.  The reader is encouraged to
consult \cite{JDH:fletcher2011} for more details about the two-dimensional case.

With exponential map, log map, and parallel transport implemented as described
here, one can perform polynomial regression on Kendall shape space in any
dimension.  We demonstrate such regression now on examples in two
dimensions. 

\subsubsection{Rat Calivarium Growth}

\begin{figure}[t]
    \begin{center}
        \includegraphics[width=1.0\textwidth]{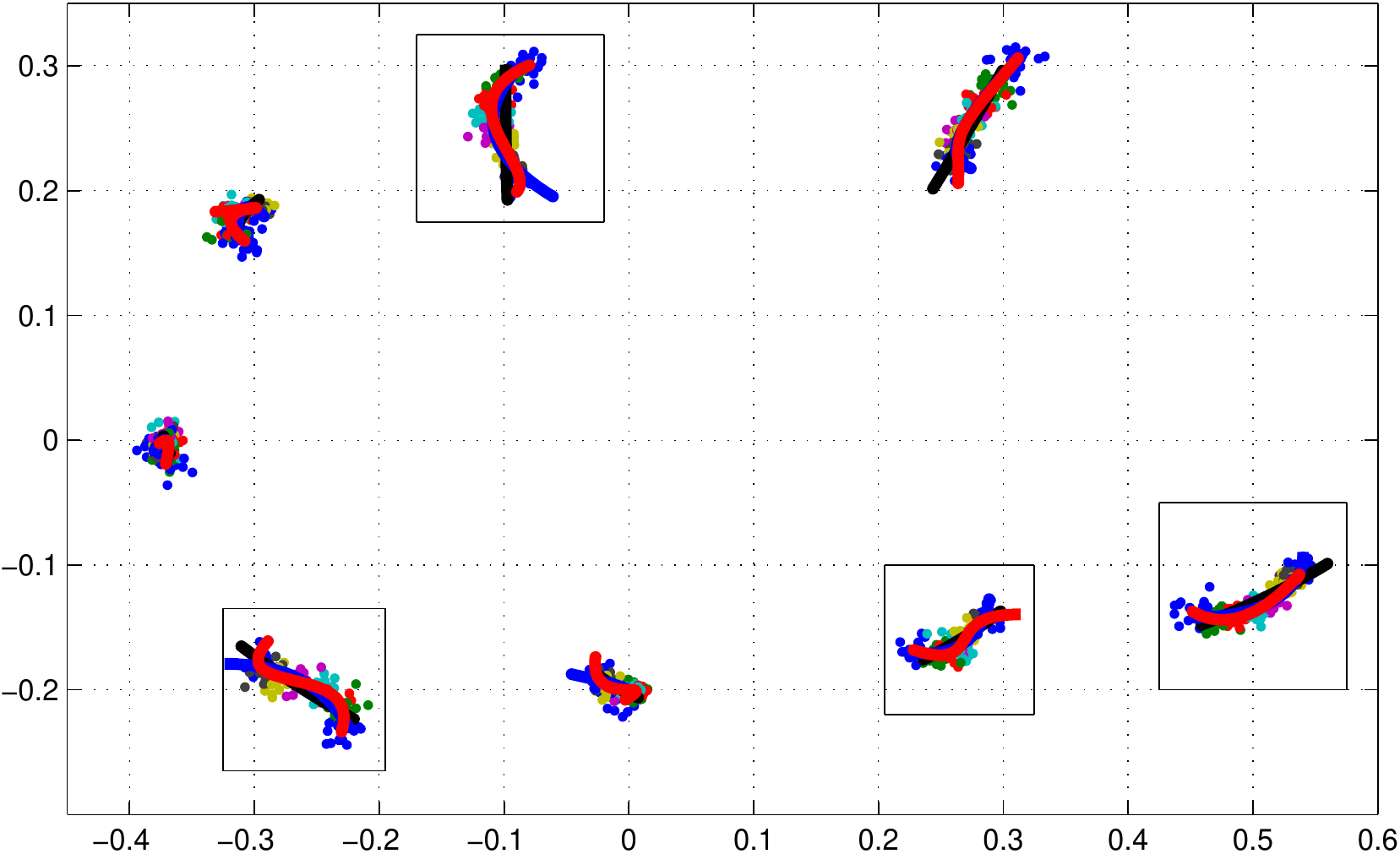}
    \end{center}
    \caption{Bookstein rat calivarium data after uniform scaling and Procrustes
    alignment.  Color indicates age group.  Zoomed views of individual
    rectangles are shown in \reffig{ratzooms}.}
    \label{fig:ratcurves}
\end{figure}

\begin{figure}[t]
    \begin{center}
        \begin{tabular}[t]{c}
            \includegraphics[width=6cm]{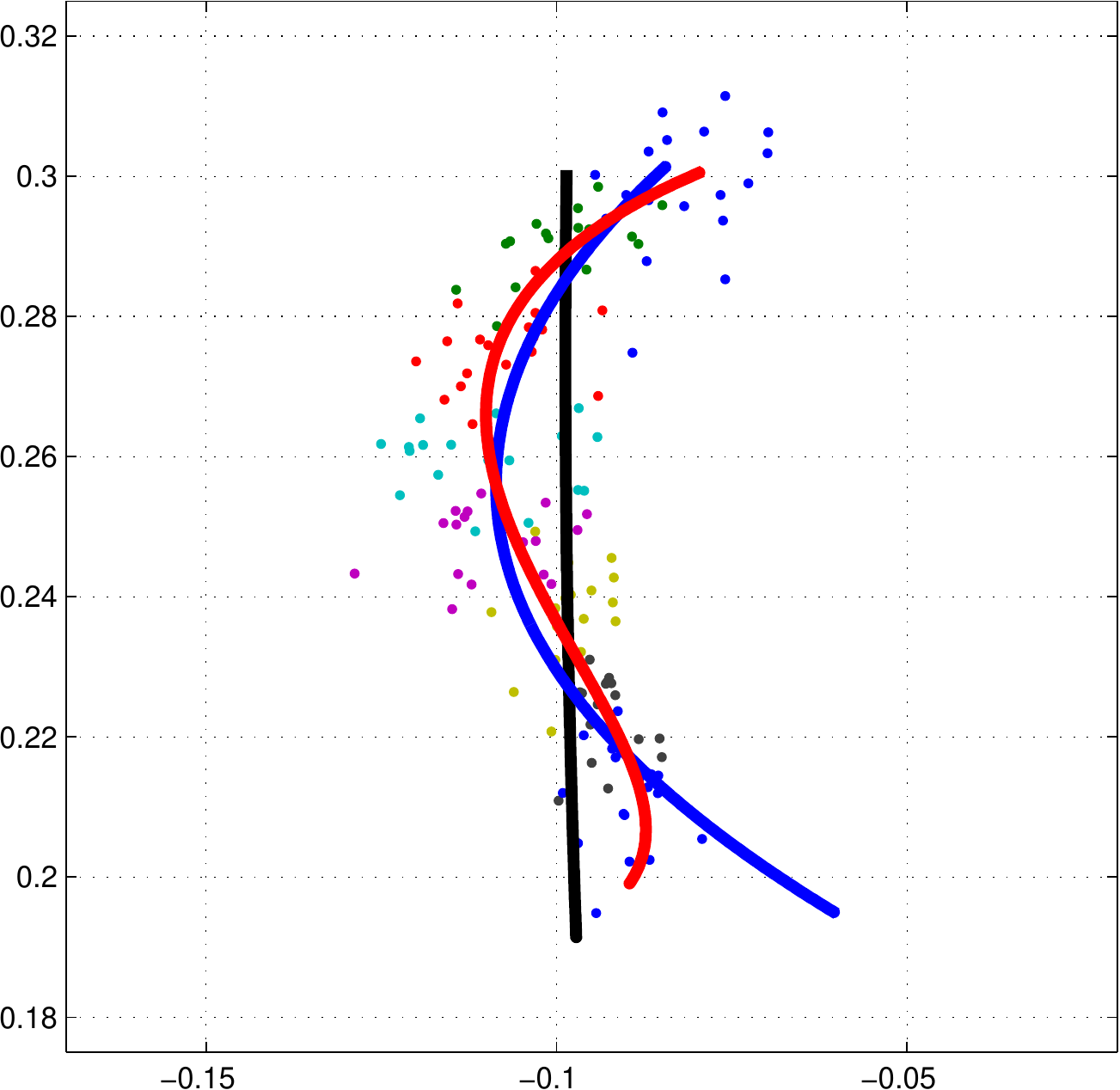} 
            \includegraphics[width=6cm]{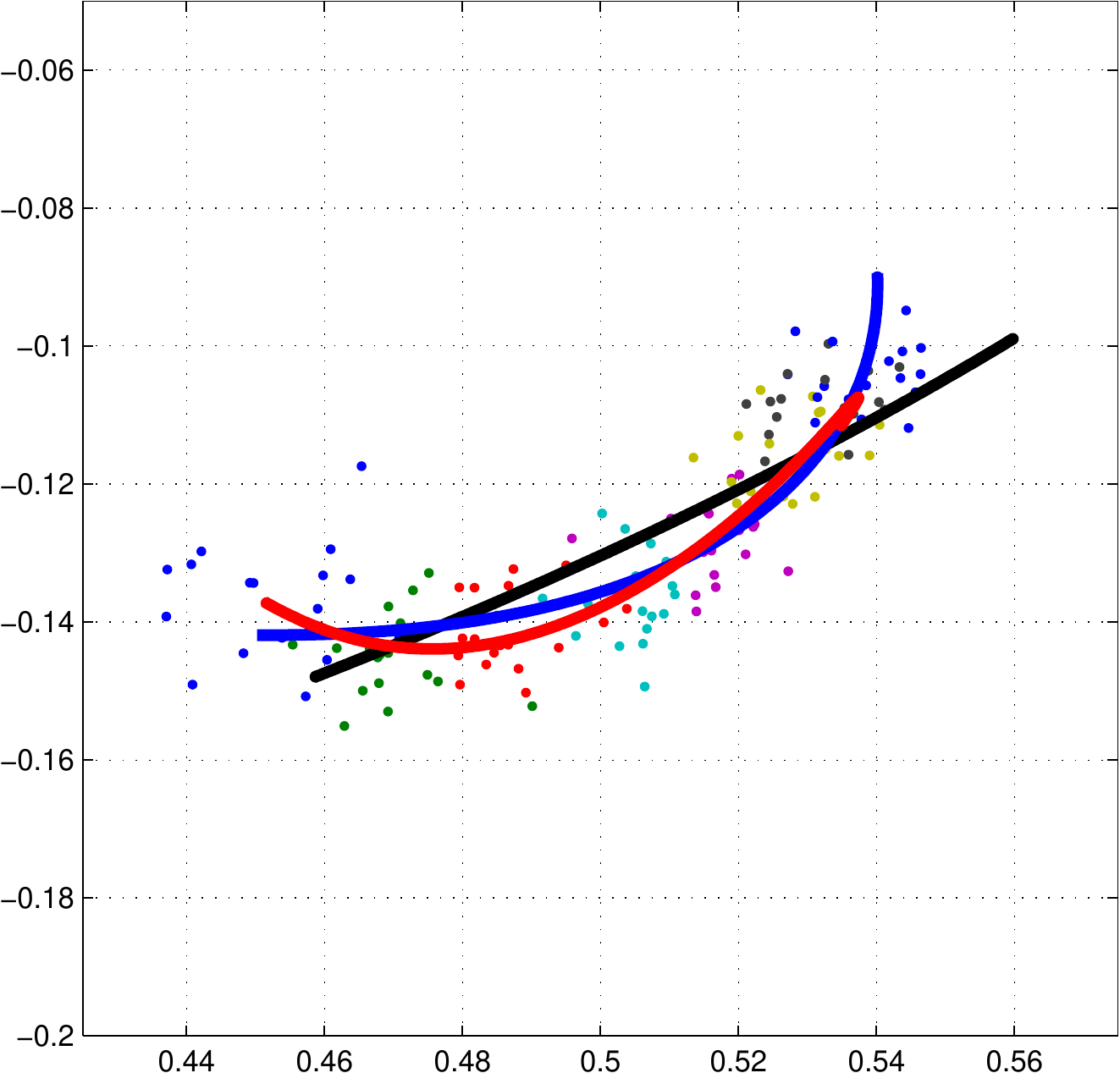} \\
            \includegraphics[width=6cm]{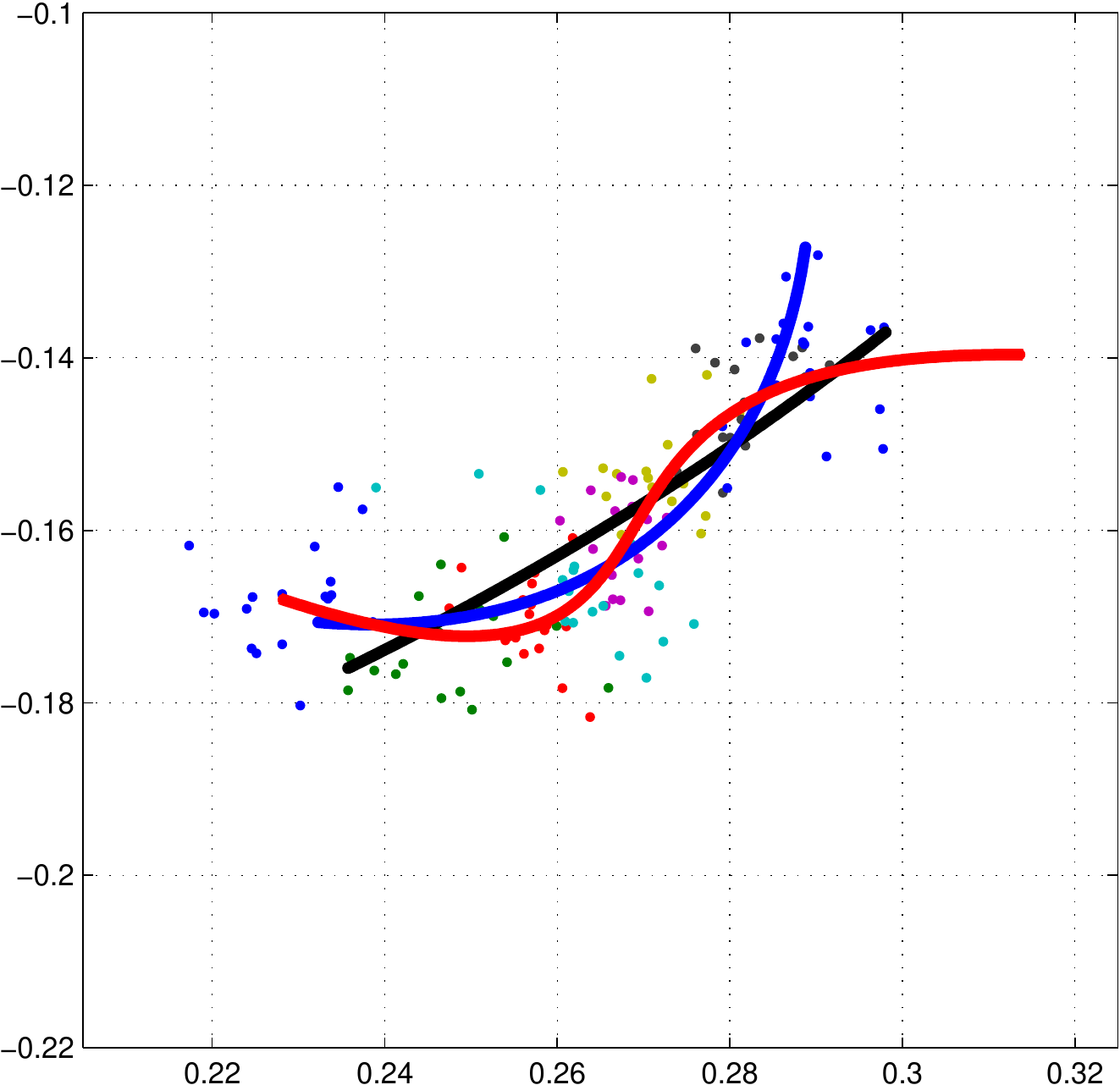} 
            \includegraphics[width=6cm]{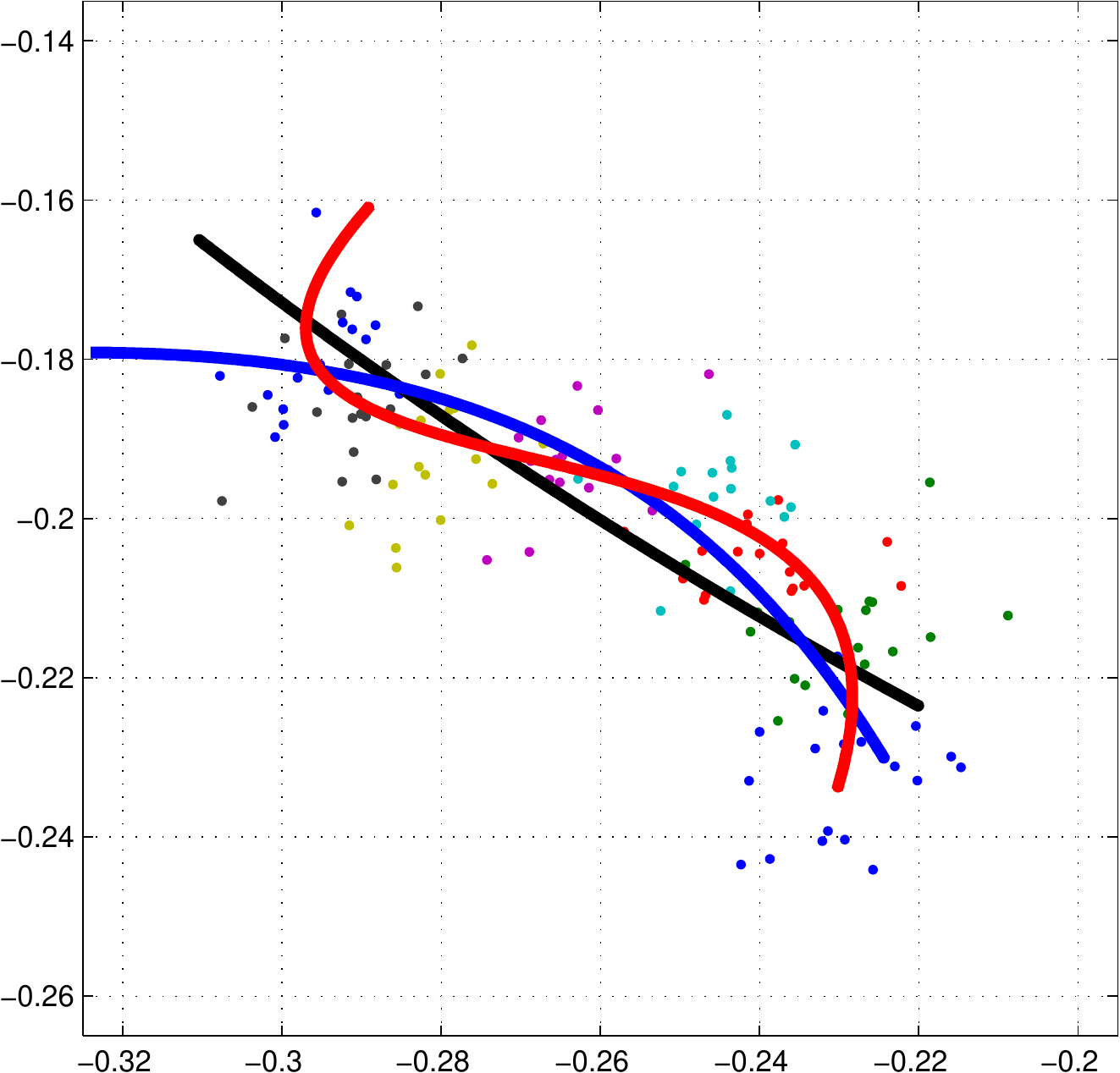}
        \end{tabular}
    \end{center}
    \caption{\label{fig:ratzooms}Zoomed views of geodesic (black, $R^2=0.79$),
    quadratic (blue, $R^2=0.85$), and cubic (red, $R^2=0.87$) regression curves
    for Bookstein rat calivarium data.}
\end{figure}

The first dataset we consider was first analyzed by
Bookstein~\cite{JDH:bookstein1991}.  The data is available for download at
\url{http://life.bio.sunysb.edu/morph/data/datasets.html} and consists of $m=8$
landmarks on a midsagittal section of rat calivaria (upper skulls).  Landmark
positions are available for 18 rats and at 8 ages apiece.  Riemannian
polynomials of orders $k=0,1,2,3$ were fit to this data.  The resulting curves
are shown in \reffig{ratcurves} and in zoomed detail in \reffig{ratzooms}.
Clearly the quadratic and cubic curves differ from that of the geodesic
regression.  The $R^2$ values agree with this qualitative difference: the
geodesic regression has $R^2=0.79$, while the quadratic and cubic regressions
both have $R^2$ values of 0.85 and 0.87, respectively.  While this shows that there is a clear improvement in the
fit due to increasing $k$ from one to two, it also shows that little is
gained by further increasing the order of the polynomial.  Qualitatively,
\reffig{ratzooms} shows
that the slight increase in $R^2$ obtained by moving from a quadratic to cubic
model corresponds to a marked difference in the curves, possibly
indicating that the cubic curve is overfitting the data.

\subsubsection{Corpus Callosum Aging}

The corpus callosum is the major white matter bundle connecting the two
hemispheres of the brain.  In order to investigate shape change of the corpus
callosum during normal aging, polynomial regression was performed on a
collection of data from the OASIS brain database (\url{www.oasis-brains.org}).
Magnetic resonance imaging (MRI) scans from 32 normal subjects with ages between 19 and
90 years were obtained from the database.  A midsagittal slice was extracted
from each volumetric image and segmented using the ITK-SNAP program
(\url{www.itksnap.org}).  Sets of 64 landmarks were optimized on the contours of these
segmentations using the ShapeWorks program~\cite{JDH:cates2007}
(\url{www.sci.utah.edu/software.html}).  The algorithm generates samplings of
each shape boundary with optimal correspondences among the populaton.



\begin{figure}[bh]
    \begin{center}
        \includegraphics[width=0.65\textwidth]{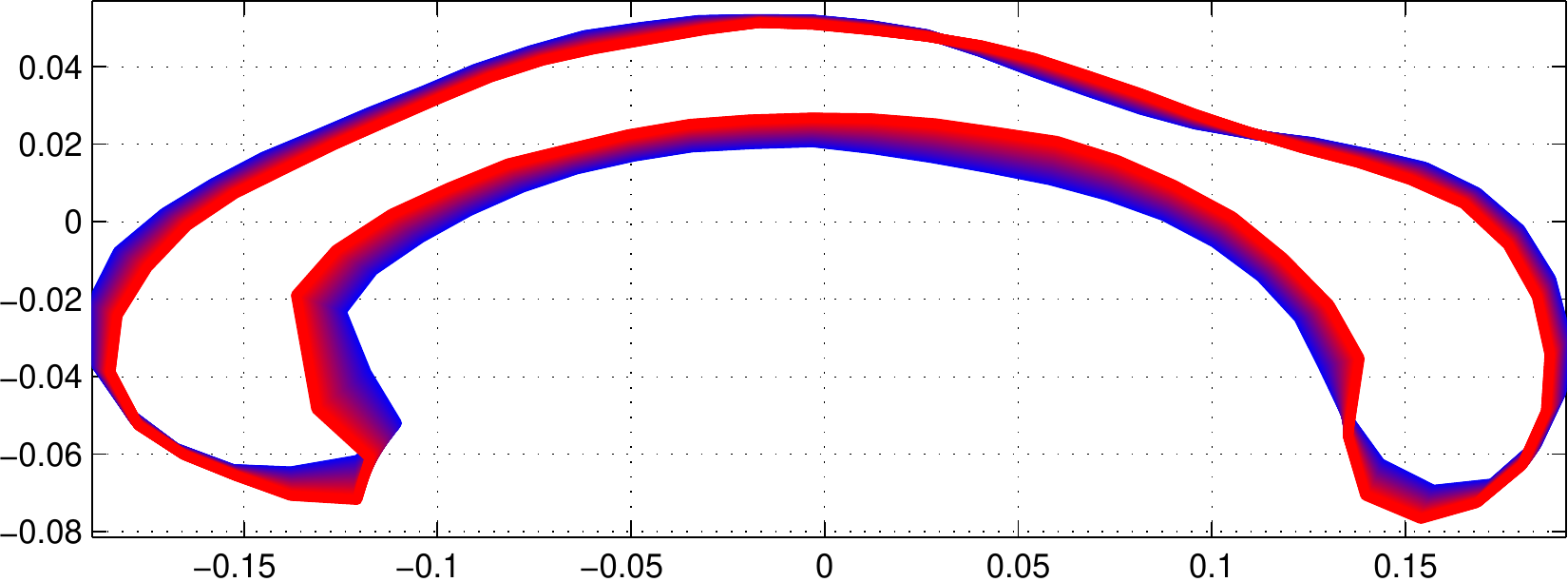} \\
        \includegraphics[width=0.65\textwidth]{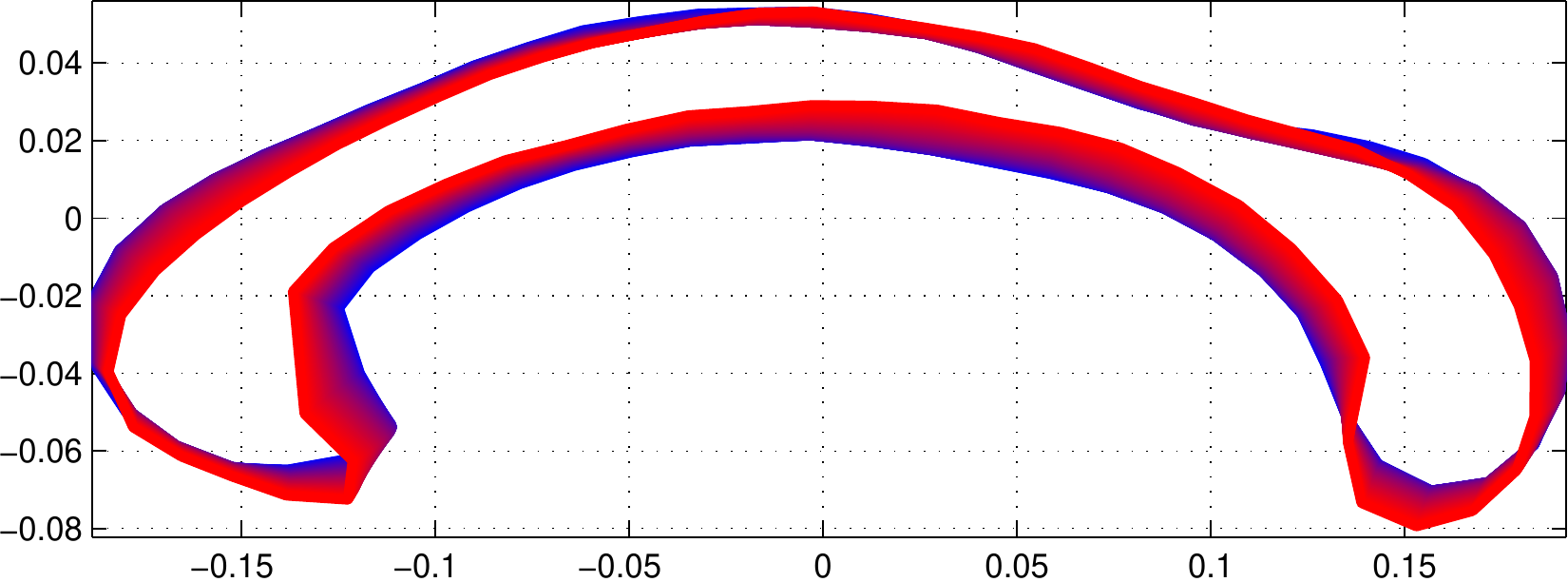} \\
        \includegraphics[width=0.65\textwidth]{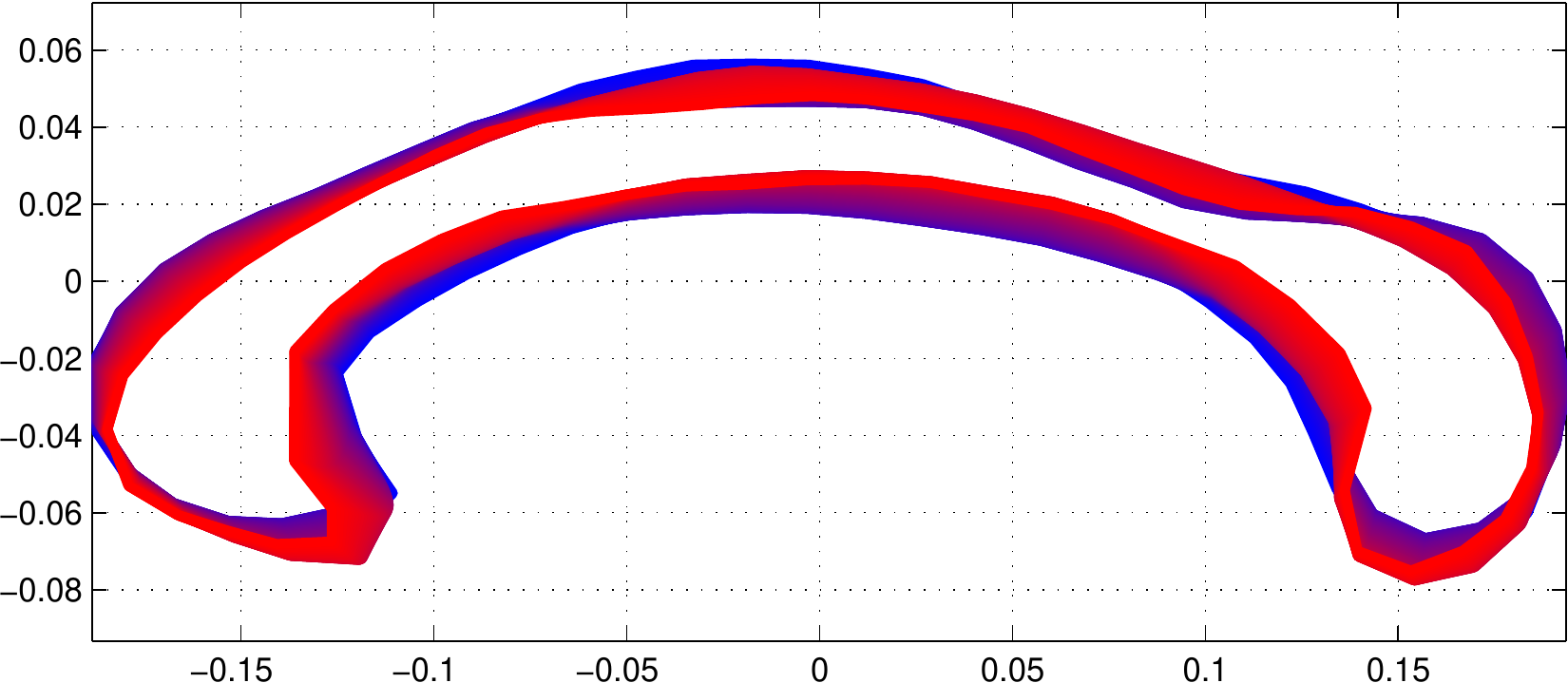}
    \end{center}
    \caption{Geodesic (top, $R^2=0.12$) quadratic (middle, $R^2=0.13$) and cubic
    (bottom, $R^2=0.21$) regression for corpus callosum dataset.  Color
    represents age, with blue indicating youth (age 19) and red indicating old
    age (age 90).}
    \label{fig:ccregressions}

\end{figure}
\begin{figure}[ht]
    \begin{center}
        \includegraphics[width=0.85\textwidth]{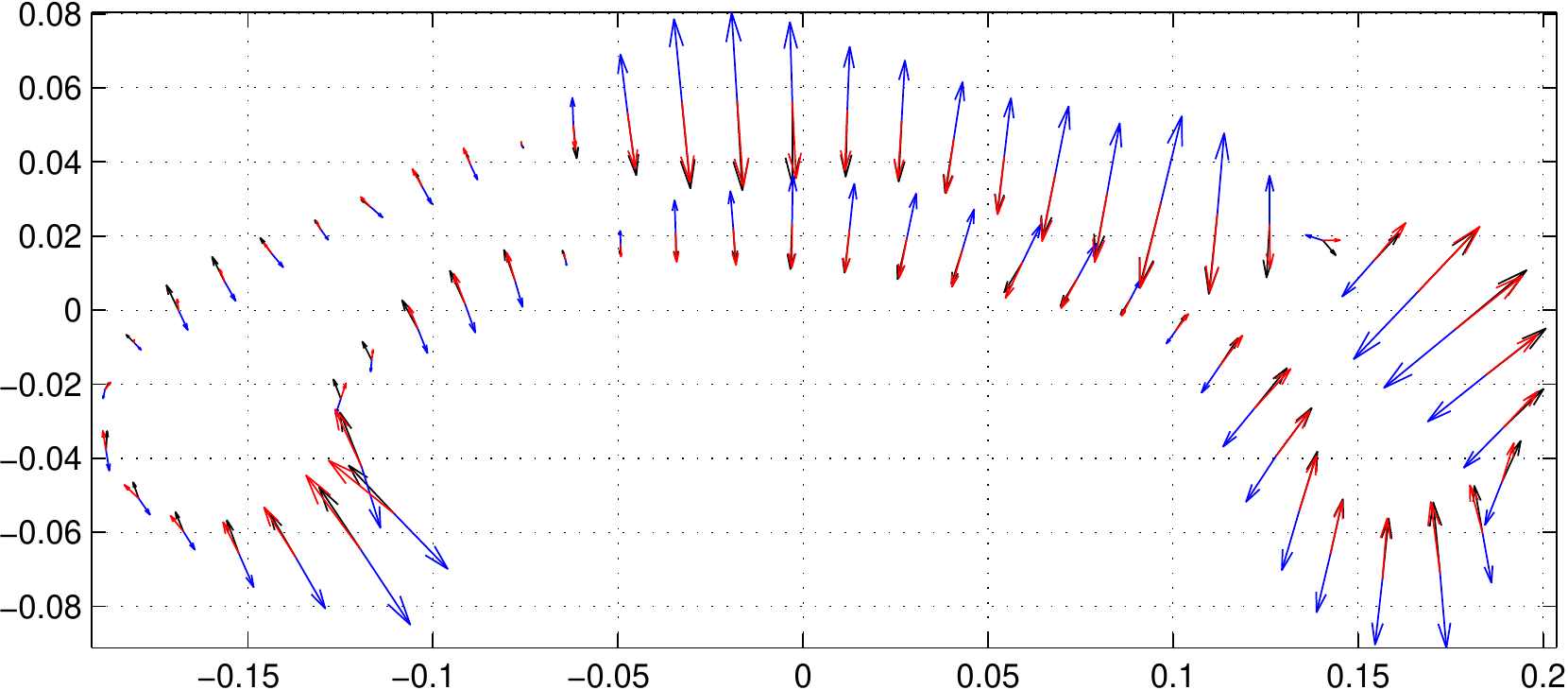}
    \end{center}
    \caption{\label{fig:ccICs}Estimated parameters for cubic regression of
    corpus callosum dataset.  The velocity (black) is nearly collinear to the
    acceleration (blue) and jerk (red).}
\end{figure}

Regression results for geodesic, quadratic, and cubic regression are shown in
\reffig{ccregressions}.  At first glance the results appear similar for the
three different models, since the motion envelopes show close agreement.
However, the $R^2$ values show an improvement from geodesic to quadratic (from
0.12 to 0.13) and from quadratic to cubic (from 0.13 to 0.21).  Inspection of
the estimated initial conditions, shown in \reffig{ccICs} reveals that the
tangent vectors appear to be rather collinear.  For the reasons stated in
\refsec{reparam}, this suggests that the differences can essentially be
described as reparametrization of the time variable, which is only accommodated
in the higher order polynomial models.  

\section*{Appendix: Derivation of Adjoint Equations}
\label{app:adj_deriv}
In this section we derive the adjoint system for the polynomial
regression problem.  The approach to calculus of variations followed
here is outlined by, for example, Noakes~\etal~\cite{JDH:noakes1989}.
Consider a simplified objective function containing only a single data
term, at time $T$:
\begin{align}
    E(\gamma,\{v_i\},\{\lambda_i\}) &= d(\gamma(T),y)^2 +\int_0^T \langle\lambda_0,\dot{\gamma}-v_1\rangle dt + \sum_{i=1}^{k-1}\int_0^T\langle
    \lambda_i,\nabla_{\dot{\gamma}}v_i-v_{i+1}\rangle dt +
    \int_0^T\langle\lambda_k,\nabla_{\dot{\gamma}}v_k\rangle dt.
\end{align}
Now consider taking variations of $E$ with respect to the vector fields
$v_i$.  For each $i$ there are only two terms containing $v_i$, so if
$W$ is a test vector field along $\gamma$, then the variation of $E$
with respect to $v_i$ in the direction $W$ satisfies 
\begin{align}
    \int_0^T\langle \delta_{v_i}E,W\rangle dt 
    &= \int_0^T\langle\lambda_i,\nabla_{\dot{\gamma}} W\rangle dt
    -\int_0^T\langle \lambda_{i-1},W\rangle dt.
\end{align}
The first term is integrated by parts to yield
\begin{align}
    \int_0^T\langle \delta_{v_i}E,W\rangle dt &=
    \langle \lambda_i, W\rangle|_0^T - \int_0^T\langle\nabla_{\dot{\gamma}} \lambda_i, W\rangle dt
    -\int_0^T\langle \lambda_{i-1},W\rangle dt.
\end{align}
The variation with respect to $v_i$ for $i=1,\ldots,k$ is then given by
\begin{align}
    \delta_{v_i(t)}E &= 0 = -\nabla_{\dot{\gamma}}\lambda_i -
    \lambda_{i-1} ,\qquad t\in(0,T) \\
    \delta_{v_i(T)}E &= 0 = \lambda_i(T) \\
    \delta_{v_i(0)}E &= -\lambda_i(t).
\end{align}
In order to determine the differential equation for $\lambda_0$, the
variation with respect to $\gamma$ must be computed.  Let $W$ again
denote a test vector field along $\gamma$.  For some $\epsilon>0$, let 
$\{\gamma_s:s\in(-\epsilon,\epsilon)\}$ be a differentiable family of
curves satisfying
\begin{align}
    \gamma_0 &= \gamma \\
    \frac{d}{ds}\gamma_s|_{s=0} &= W.
\end{align}
If $\epsilon$ is chosen small enough, the vector field $W$ can be
extended to a neighborhood of $\gamma$ such that
$[W,\dot{\gamma_s}]=0$, where a dot indicates the derivative in the
$\frac{\partial}{\partial t}$ direction.  The vanishing Lie bracket
implies the following identities
\begin{align}
    \nabla_W\dot{\gamma}_s &= \nabla_{\dot{\gamma}_s}W \\
    \nabla_W\nabla_{\dot{\gamma}_s} &= \nabla_{\dot{\gamma}_s}\nabla_W +
    R(W,\dot{\gamma}_s).
\end{align}
Finally, the vector fields $v_i,\lambda_i$ are extended along
$\gamma_s$ via parallel translation, so that
\begin{align}
    \nabla_W v_i &= 0 \\
    \nabla_W \lambda_i &= 0.
\end{align}
The variation of $E$ with respect to $\gamma$ satisfies
\begin{align}
    \int_0^T\langle \delta_\gamma E,W\rangle dt &=
    \frac{d}{ds}E(\gamma_s,\{v_i\},\{\lambda_i\})|_{s=0} \\
    &= -\langle\Log_{\gamma(T)}y,W(T)\rangle + \frac{d}{ds}\int_0^T \langle\lambda_0,\dot{\gamma_s}-v_1\rangle dt|_{s=0} \\
    &\quad + \frac{d}{ds}\sum_{i=1}^{k-1}\int_0^T\langle
    \lambda_i,\nabla_{\dot{\gamma}_s}v_i-v_{i+1}\rangle dt|_{s=0} \\
    &\quad+
    \frac{d}{ds}\int_0^T\langle\lambda_k,\nabla_{\dot{\gamma}_s}v_k\rangle dt|_{s=0}.
\end{align}
As the $\lambda_i$ are extended via parallel translation, their inner
products satisfy
\begin{align}
    \frac{d}{ds}\langle \lambda_i, U\rangle|_{s=0} &= \langle \nabla_W\lambda_i,U\rangle + \langle
    \lambda_i,\nabla_W U\rangle = \langle \lambda_i,\nabla_W U\rangle.
\end{align}
Then applying this to each term in the previous equation,
\begin{align}
    \int_0^T\langle \delta_\gamma E,W\rangle dt &= -\langle\Log_{\gamma(T)}y,W(T)\rangle + \int_0^T \langle\lambda_0,\nabla_W\dot{\gamma}-\nabla_Wv_1\rangle dt \\
    &\quad + \sum_{i=1}^{k-1}\int_0^T\langle
    \lambda_i,\nabla_W\nabla_{\dot{\gamma}}v_i-\nabla_Wv_{i+1}\rangle dt \\
    &\quad+
    \int_0^T\langle\lambda_k,\nabla_W\nabla_{\dot{\gamma}}v_k\rangle dt.
\end{align}
Then by construction, since $\nabla_W v_i=0$,
\begin{align}
    \int_0^T\langle \delta_\gamma E,W\rangle dt &= -\langle\Log_{\gamma(T)}y,W(T)\rangle + \int_0^T \langle\lambda_0,\nabla_W\dot{\gamma}\rangle dt
     + \sum_{i=1}^{k}\int_0^T\langle
    \lambda_i,\nabla_W\nabla_{\dot{\gamma}}v_i\rangle dt.
\end{align}
Then using the Lie bracket and curvature identities, this is written as
\begin{align}
    \int_0^T\langle \delta_\gamma E,W\rangle dt &= -\langle\Log_{\gamma(T)}y,W(T)\rangle + \int_0^T \langle\lambda_0,\nabla_{\dot{\gamma}}W\rangle dt
    + \sum_{i=1}^{k}\int_0^T\langle
    \lambda_i,\nabla_{\dot{\gamma}}\nabla_W v_i + R(W,\dot{\gamma})v_i\rangle dt,
\end{align}
which is further simplified, again using the identity $\nabla_W v_i=0$:
\begin{align}
    \int_0^T\langle \delta_\gamma E,W\rangle dt &= -\langle\Log_{\gamma(T)}y,W(T)\rangle + \int_0^T \langle\lambda_0,\nabla_{\dot{\gamma}}W\rangle dt
    + \sum_{i=1}^{k}\int_0^T\langle
    \lambda_i,R(W,\dot{\gamma})v_i\rangle dt,
\end{align}
Using the Bianchi identities, it can be demonstrated that the
curvature tensor satisfies the identity~\cite{JDH:docarmo1992}:
\begin{align}
    \langle A,R(B,C)D\rangle &= -\langle B,R(D,A)C\rangle,
\end{align}
for any vectors $A,B,C,D$.  The covariant derivative along $\gamma$ is
also integrated by parts to arrive at
\begin{align}
    \int_0^T\langle \delta_\gamma E,W\rangle dt &= -\langle\Log_{\gamma(T)}y,W(T)\rangle + \langle\lambda_0,W\rangle|_0^T - \int_0^T \langle\nabla_{\dot{\gamma}}\lambda_0,W\rangle dt
    - \sum_{i=1}^{k}\int_0^T\langle
    R(v_i,\lambda_i)\dot{\gamma},W\rangle dt.
\end{align}
Finally, gathering terms, the adjoint equation for $\lambda_0$ and its
gradients are obtained:
\begin{align}
    \delta_{\gamma(t)}E &= 0 = -\nabla_{\dot{\gamma}}\lambda_0 -
    \sum_{i=1}^k R(v_i,\lambda_i)\dot{\gamma},\qquad t\in(0,T) \\
    \delta_{\gamma(T)}E &= 0 = -\Log_{\gamma(T)}y + \lambda_0 \\
    \delta_{\gamma(0)}E &= -\lambda_0.
\end{align}
Along with the variations with respect to $v_i$, this constitutes the
full adjoint system.  Extension to the case of multiple data at multiple
time points is trivial, and results in the adjoint system presented in
\refsec{estimation}.



{\small
\bibliographystyle{apalike}
\bibliography{arxiv}
}

\end{document}